\documentclass{article}

\usepackage{amssymb}

\newtheorem{lem}{Lemma}

\newtheorem{prop}{Proposition}

\begin{document}

\title{The Roman harmonic numbers revisited}

\author{J. Sesma\thanks{Email: javier@unizar.es} \\ \ \\
{\em Departamento de F\'{\i}sica Te\'{o}rica, Facultad de Ciencias,} \\
{\em 50009, Zaragoza, Spain}}

\maketitle

\begin{abstract}
Two decades ago, Steven Roman, Daniel E. Loeb and Gian-Carlo Rota introduced a family of harmonic numbers in their study of harmonic logarithms. We propose to refer to those numbers as {\it Roman harmonic numbers}. With the purpose of revitalizing the study of these mathematical objects, we recall here their known properties and unveil additional ones. An integral representation, several generating relations, and a collection of sum rules involving those numbers are presented. It is also shown that higher derivatives of the Pochhammer and reciprocal Pochhammer symbols are easily expressed in terms of Roman harmonic numbers.
\end{abstract}



\section{Introduction}

In the course of a research on the so called epsilon expansions of Appell and Kamp\'e de F\'eriet functions \cite{grey}, we considered convenient to define what we called modified generalized harmonic numbers, $\hat{H}_n^{(k)}$, in the form
\begin{equation}
\hat{H}_0^{(k)}=\delta_{k,0}\,, \qquad  \hat{H}_n^{(k)} \equiv \sum_{j=1}^n (-1)^{j-1}\,{n \choose j}\,\frac{1}{j^k}\,, \quad n\geq 1\,.  \label{I1}
\end{equation}
Later on, we realized that these mathematical objects, with the name of harmonic numbers and denoted by $c_n^{(k)}$, had already been introduced, by Loeb and Rota \cite{loeb} and by Roman \cite{rom1,rom2}, to give explicit expressions of the harmonic logarithms. On the other hand, the moments
\begin{equation}
\bar{d}^p = \int_0^\infty d\delta\,\delta^p\,\mathcal{P}(\delta)\,  \qquad p=0, 1, 2, \ldots\,,   \label{C1}
\end{equation}
of the quantum probability distribution
\begin{equation}
\mathcal{P}(\delta) = 2\left(D_\varepsilon-1\right)\,e^{-2\delta}\left(1-e^{-2\delta}\right)^{D_\varepsilon-2}\,,  \quad D_\varepsilon =2, 3, \dots\,, \label{C2}
\end{equation}
discussed by Coffey \cite{cof1} and used (for $p=1$ and $2$) to approximately calculate the decoherence factors for $D_\varepsilon$-dimensional quantum states, are closely related to the same harmonic numbers, as shown below.
In spite of their interesting properties, the $c_n^{(k)}$ have not received, in our opinion, due attention by number theorists. Here we intend to arouse the interest in those harmonic numbers by recalling their known properties and presenting some others.

Several generalizations of harmonic numbers have been used by different authors \cite{cand,cho2,cho3,cho4,cho1,copp,kron,liuw,sofo}.
To avoid misunderstandings, we adopt the name of {\em Roman harmonic numbers} for those introduced in Refs.~\cite{loeb}, \cite{rom1}, and \cite{rom2}.
One of the possible equivalent definitions of these numbers is given in Section 2, which recalls also their main properties, found by Roman \cite{rom1,rom2} and by Loeb and Rota \cite{loeb}. Section 3 shows, firstly, the connection of those numbers with known nested sums. Their integral representation allows to discover the relation between Roman harmonic numbers and Coffey's quantum distribution moments. Then, several generating relations, sum rules, and additional properties not mentioned before are presented. As an application, expressions are given of the derivatives of the Pochhammer and reciprocal Pochhammer symbols in terms of the $c_n^{(k)}$. To end, some pertinent comments are included in Section 4.

\section{First definitions and known properties}

Before recalling  the Roman harmonic numbers, some auxiliary definitions, given in \cite{loeb}, \cite{rom1}, and \cite{rom2}, are necessary. We present them with the notation used in these references.

\noindent{\bf Definition 1.} (Roman number.)
The {\em Roman} $n$ is defined to be
\begin{equation}
\lfloor n\rceil  = \left\{ \begin{array}{lll}n & {\rm for} & n\neq 0\,, \\ 1 & {\rm for} & n=0\,. \end{array} \right.  \label{II1}
\end{equation}

\noindent{\bf Definition 2.} (Roman factorial.)
For every integer $n$, define $n$ {\em Roman factorial} to be
\begin{equation}
\lfloor n\rceil ! = \left\{ \begin{array}{lll}(-1)^{n+1}/(-n-1)! & {\rm for} & n<0\,, \\ n! & {\rm for} & n\geq 0\,. \end{array} \right.  \label{II2}
\end{equation}

From definitions 1 and 2, it is immediate to check that \cite[Prop. 3.2]{rom2}
\begin{equation}
\lfloor n\rceil ! = \lfloor n\rceil \,  \lfloor n-1\rceil !\,.  \label{II3}
\end{equation}

As we are going to see, the Roman harmonic numbers are closely related to the {\em Stirling numbers of the first kind}, $s(n,k)$. These are commonly defined for integer $n\geq 0$ by means of their ordinary generating function \cite{rio1,rio2,temm}. A different, although equivalent, definition due to Loeb and Rota \cite[Def. 3.3.2]{loeb} allows to consider also negative values of $n$. Moreover, Coffey \cite[Appendix B]{coff} has given a definition of generalized Stirling numbers with complex first argument, which is consistent with that proposed by Loeb and Rota. Nevertheless, along this paper we restrict ourselves to the case of integer first argument.

\noindent{\bf Definition 3.} (Stirling numbers of the first kind.) For all integers $n$ and non-negative integers $k$, the {\em Stirling numbers} $s(n,k)$ are defined as the coefficients of the formal expansion
\begin{equation}
\frac{\Gamma(x+1)}{\Gamma(x+1-n)}=\sum_{k=0}^\infty\,s(n,k)\,x^k\,.   \label{II4}
\end{equation}

\noindent{\bf Remark 1.}
The left-hand side of (\ref{II4}) becomes more concise when written in terms of Pochhammer symbols. Different conventions are used to refer to such objects. Here we adopt that of the Digital Library of Mathematical Functions \cite[Eq.~5.2.4]{dlmf}, namely
\[
(z)_0=1\,,   \qquad  (z)_n=z(z+1)\cdots(z+n-1)\,.
\]
With this convention, Eq.~(\ref{II4}) splits into
\begin{equation}
(x)_n=(-1)^n\,\sum_{k=0}^\infty\,s(n,k)\,(-x)^k\,,  \qquad \frac{1}{(x+1)_n}=\sum_{k=0}^\infty\,s(-n,k)\,x^k\,, \qquad n\geq 0\,.
    \label{II5}
\end{equation}
Since $(x)_n$ is a polynomial of degree $n$ in $x$, $s(n,k)$ vanishes for $k>n\geq 0$. Instead, the $s(-n,k)$, with $n>0$, may be different from zero for arbitrarily large values of $k$.
From the obvious expression
\[
s(n,k)=\frac{1}{k!}\,\left.\frac{d^k}{dx^k}\left(\frac{\Gamma(x+1)}{\Gamma(x+1-n)}\right)\right|_{x=0}\,,
\]
it follows immediately
\begin{equation}
s(0,k)=\delta_{k,0}\,,  \qquad s(n,0)=\left\{ \begin{array}{lll} 1/(-n)! & {\rm for} & n<0 \,,  \\ 0 & {\rm for} & n>0\,. \end{array} \right.     \label{II7}
\end{equation}

\begin{lem}\hspace{-4pt}{\bf .}
For all integer $n$, and positive integer $k$, the Stirling numbers of the first kind obey the recurrence relation
\begin{equation}
n\,s(-n,k)=-\,s(-n,k-1)+s(-n+1,k)\,.  \label{II8}
\end{equation}
\end{lem}
\noindent \textsc{Proof.}
Compare coefficients of equal powers of $x$ in the expression obtained by expanding, according to (\ref{II4}), the three fractions in the evident relation
\begin{equation}
n\,\frac{\Gamma(x+1)}{\Gamma(x+1+n)} = -\,x\,\frac{\Gamma(x+1)}{\Gamma(x+1+n)} + \frac{\Gamma(x+1)}{\Gamma(x+n)}\,. \label{II9}
\end{equation}
\hspace{330pt} $\square$

Now we recall the definition of the $c_n^{(k)}$, as given in \cite[Prop. 6.1]{rom2}.

\noindent{\bf Definition 4.} (Roman harmonic numbers).
For all integers $n$ and nonnegative integers $k$, the {\em Roman harmonic numbers} $c_n^{(k)}$, of order $k$ and degree $n$, are uniquely determined by the initial conditions
\begin{equation}
 c_0^{(k)}= \delta_{k,0}\,,  \qquad c_n^{(0)}=\left\{ \begin{array}{lll} 1 & {\rm for} & n\geq 0 \,,  \\ 0 & {\rm for} & n<0\,, \end{array} \right.   \label{II10}
\end{equation}
and the recurrence relation (for $k>0$)
\begin{equation}
n\,c_n^{(k)}=c_n^{(k-1)}+\lfloor n\rceil\,c_{n-1}^{(k)}\,,  \label{II11}
\end{equation}
used for increasing absolute values of $n$.

From this definition, it is easy to check that
\begin{equation}
c_1^{(k)}=1\, \qquad  {\rm and}  \qquad  c_{-n}^{(1)}=-1 \quad {\rm for}  \quad n>0\,.    \label{II12}
\end{equation}
Besides,
\[
c_{-n}^{(k)}\leq 0 \qquad {\rm and} \qquad c_n^{(k)}>0 \quad {\rm for} \quad n>0\,.
\]
Numerical values of the $c_n^{(k)}$  for $-6\leq n\leq 5$ and $0\leq k\leq 6$ are given in \cite[Table 3.9]{loeb}.

The relation of the Roman harmonic  to the Stirling numbers is given by the following proposition, which states a property which has been used in \cite[Def. 3.3.3]{loeb} to define the $c_n^{(k)}$.

\begin{prop}\hspace{-4pt}{\bf .}
For all integers $n$ and all non-negative integers $k$,
\begin{equation}
c_n^{(k)}=(-1)^k\,\lfloor n\rceil !\, s(-n,k)\,.       \label{II13}
\end{equation}
\end{prop}

\noindent \textsc{Proof.}
For $n=0$, compare Eqs.~(\ref{II7}) and (\ref{II10}). For $n\neq 0$, (\ref{II13}) follows from comparison of (\ref{II8}) and (\ref{II11}), bearing in mind (\ref{II3}). \hspace{112pt} $\square$

The relation (\ref{II13}) allows to write a large list of properties of the $c_n^{(k)}$ with nonpositive $n$ as immediate translations of the widely studied properties of the Stirling numbers of the first kind with nonnegative first argument. We do not reproduce here the resulting relations, for the sake of brevity, but refer the reader to, for instance, \cite[Chap.~V]{comt}, or \cite[\S 26.8]{dlmf}. Nevertheless, we report here the following ones
\begin{equation}
c_{-n}^{(k)}=0 \qquad {\rm if}\quad  k>n\geq 0,    \label{II14}
\end{equation}
\begin{equation}
{k \choose h}\frac{c_{-n}^{(k)}}{-n} = \sum_{j=k-h}^{n-h}\frac{c_{-n+j}^{(h)}}{-n+j}\frac{c_{-j}^{(k-h)}}{-j}\,,
  \qquad n \geq k > h \geq 0\,,  \label{II15}
\end{equation}
\begin{equation}
c_{-n-1}^{(k+1)}=\sum_{j=k}^n\frac{1}{j}\,c_{-j}^{(k)}\,,  \qquad  n \geq k > 0\,,   \label{II16}
\end{equation}
to be used below.

Properties of the $c_n^{(k)}$ with positive $n$ are, instead, no so trivially obtained. Section 3 of Ref.~\cite{loeb} contains several explicit expressions of the $c_n^{(k)}$ and some of their properties. As stated before, our purpose is to
recall those that are known and to point out some others which we have encountered. In what follows, we focus on the Roman harmonic numbers $c_n^{(k)}$ of positive degree $n$.

The fact that our before mentioned modified  generalized harmonic numbers $\hat{H}_n^{(k)}$ are but the Roman harmonic numbers $c_n^{(k)}$, with nonnegative $n$, becomes evident from comparison of Eqs.~(\ref{I1}) and (\ref{II10}), and from the second of the properties mentioned in the following

\begin{prop}\hspace{-4pt}{\bf .} {\rm (\cite[Prop. 6.3]{rom2}.)}
The Roman harmonic numbers  $c_n^{(k)}$  with $n>0$ have the following properties:
\begin{description}
\item [{\rm (a)}]
\begin{equation}
c_n^{(k)}=\sum_{j=1}^n\ \frac{1}{j}\,c_j^{(k-1)}\,,  \qquad k\geq 1\,. \label{II17}
\end{equation}
\item [{\rm (b)}]
\begin{equation}
c_n^{(k)}=\sum_{j=1}^n (-1)^{j-1}\,{n \choose j}\,\frac{1}{j^k}\,.  \label{II18}
\end{equation}
\item [{\rm (c)}] For each $n>1$, the sequence $c_n^{(k)}$ forms a strictly increasing sequence in $k$ such that
\begin{equation}
\lim_{k\to\infty}c_n^{(k)}=n\,.      \label{II19}
\end{equation}
\end{description}
\end{prop}

\noindent{\bf Remark 2.}
It follows from (\ref{II10}) and (\ref{II17}) that the $c_n^{(1)}$ are the familiar harmonic numbers,
\begin{equation}
c_n^{(1)} = \sum_{j=1}^n \frac{1}{j} = H_n\,.    \label{II20}
\end{equation}

\section{Additional properties}

\subsection{Relation to nested sums}

The property (a) in Proposition 2 shows that the $c_n^{(k)}$ with $n\geq 1$ are particular cases of the nested $S$-sums \cite{moch} defined by
\[
S(n)  =  \left\{ \begin{array}{lll} 0 & {\rm for} & n\leq 0 \,,  \\ 1 & {\rm for} & n>0\,,\end{array}\right.
\]
\[
S(n;m_1,\dots,m_k;x_1,\ldots,x_k)  =  \sum_{j=1}^n\frac{x_1^j}{j^{m_1}}\,S(j;m_2,\ldots,m_k;x_2,\ldots,x_k)\,.
\]
In fact,
\begin{equation}
c_n^{(k)}=S(n;1,\ldots_{(k)}\ldots,1;1,\ldots_{(k)}\ldots,1)\,,   \qquad  n\geq 1\,.  \label{IIIi1}
\end{equation}
On the other hand, for the Roman harmonic numbers of negative degree, comparison of Eqs.~(\ref{II12}) and (\ref{II16}) with the definition of the nested $Z$-sums \cite{moch}
\[
Z(n)  =  \left\{ \begin{array}{lll} 0 & {\rm for} & n<0 \,,  \\ 1 & {\rm for} & n\geq 0\,,\end{array}\right.
\]
\[
Z(n;m_1,\dots,m_k;x_1,\ldots,x_k)  =  \sum_{j=1}^n\frac{x_1^j}{j^{m_1}}\,Z(j-1;m_2,\ldots,m_k;x_2,\ldots,x_k)\,,
\]
shows that
\begin{equation}
c_{-n}^{(k)}=-\,Z(n-1;1,\ldots_{(k-1)}\ldots,1;1,\ldots_{(k-1)}\ldots,1)\,,   \qquad  n\geq 1\, \quad k\geq 1\,.  \label{IIIi3}
\end{equation}
In the notation used by Vermaseren \cite{verm} for the harmonic sums and the Euler-Zagier sums, the relations (\ref{IIIi1}) and (\ref{IIIi3}) become, respectively,
\begin{eqnarray}
c_n^{(k)} & = & S_{1,\ldots_{(k)}\ldots,1}(n)\,,  \qquad n\geq 1\,,   \label{IIIi4} \\
c_{-n}^{(k)} & = & -\,Z_{1,\ldots_{(k-1)}\ldots,1}(n-1)\,,   \qquad  n\geq 1\,, \quad k\geq 1\,.  \label{IIIi5}
\end{eqnarray}

Vermaseren has also pointed out that the harmonic sum $S_{1,\ldots_{(k)}\ldots,1}(n)$, and therefore $c_n^{(k)}$, can be written as a sum of products of generalized harmonic numbers $H_n^{(k)}$ defined by
\begin{equation}
H_n^{(k)} = \sum_{j=1}^n\frac{1}{j^k}\,.  \nonumber
\end{equation}
With this purpose, let us consider all possible different decompositions
of $k$ as a sum of products of two positive integers, for instance
\begin{equation}
d_i: \qquad k=m_{1,i}\,l_{1,i} + m_{2,i}\,l_{2,i} + \ldots + m_{n_i,i}\,l_{n_i,i}\,,  \nonumber
\end{equation}
where the $l_{j,i}$ ($j=1,2,\ldots, n_i$) are different.
Each one of these decompositions may be associated to the conjugation class of the symmetric group $S_k$ constituted by the elements characterized by containing $m_{1,i}$ cycles of length $l_{1,i}$, $m_{2,i}$ cycles of length $l_{2,i}$, $\ldots$, and $m_{n_i,i}$ cycles of length $l_{n_i,i}$. As it is well known, the number of elements of such conjugation class is
\[
k!\,\prod_{j=1}^{n_i}\,\frac{1}{m_{j,i}!\,l_{j,i}^{m_{j,i}}}
\]Then,
\begin{equation}
c_n^{(k)} = \sum_i\,\left( \prod_{j=1}^{n_i}\,\frac{1}{m_{j,i}!}\left(\frac{H_n^{(l_{j,i})}}{l_{j,i}}\right)^{m_{j,i}}\right)\,,
  \label{IIIi6}
\end{equation}
where the sum extends to all different decompositions $d_i$. As an example,
\begin{eqnarray}
c_n^{(5)} & = & \frac{1}{120}\,\left(H_n^{(1)}\right)^5 + \frac{1}{12}\,\left(H_n^{(1)}\right)^3 H_n^{(2)} +  \frac{1}{6}\,\left(H_n^{(1)}\right)^2H_n^{(3)} \nonumber  \\
& & +\ \frac{1}{8}\,H_n^{1}\left(H_n^{(2)}\right)^2 + \frac{1}{4}\,H_n^{(1)}H_n^{(4)}  + \frac{1}{6}\,H_n^{(2)}H_n^{(3)} + \frac{1}{5}\,H_n^{(5)}\,.     \nonumber
\end{eqnarray}

\subsection{Integral representation}

A known integral representation of the harmonic sums \cite[Eq.~(35)]{verm} is valid, in view of (\ref{IIIi4}), for the Roman harmonic numbers, as the following lemma states.

\begin{lem}\hspace{-4pt}{\bf .} {\rm (Integral representation).}
The Roman harmonic numbers of positive degree admit the integral representation
\begin{equation}
c_n^{(k)} = (-1)^k\,\frac{n}{k!}\int_0^1\,dx\,x^{n-1}\,[\ln (1-x)]^k\,,  \qquad n\geq 1\,.  \label{IIIi7}
\end{equation}
\end{lem}

\noindent \textsc{Proof.}
A trivial change of integration variable gives for the integral in the right-hand side of (\ref{IIIi7})
\begin{eqnarray}
\int_0^1\,dx\,x^{n-1}\,[\ln (1-x)]^k & = & \int_0^1\,dy\,(1-y)^{n-1}\,[\ln y]^k     \nonumber   \\
   & = & \sum_{j=0}^{n-1}{n-1 \choose j}\int_0^1\,dy\,(-y)^{j}\,[\ln y]^k\,.   \label{IIIi8}
\end{eqnarray}
The value of the last integral can be found by repeated integration by parts or taken from \cite[Eq.~2.6.3.2]{prud}. One obtains in this way
\begin{equation}
\int_0^1\,dx\,x^{n-1}\,[\ln (1-x)]^k = \sum_{j=0}^{n-1}{n-1 \choose j}(-1)^k\,k!\,\frac{(-1)^j}{(j+1)^{k+1}}\,,
\qquad n\geq 1\,,     \label{IIIi9}
\end{equation}
which can be written in the form
\begin{equation}
\int_0^1\,dx\,x^{n-1}\,[\ln (1-x)]^k = \frac{(-1)^k\,k!}{n}\sum_{l=1}^{n}{n \choose l}(-1)^{l-1}\,\frac{1}{l^{k}}\,,
\qquad n\geq 1\,,     \label{IIIi10}
\end{equation}
Comparison of this equation with (\ref{II18}) proves the lemma. \hspace{80pt}$\square$

\noindent{\bf Remark 3.}
It is evident, from the recurrence relation (\ref{II11}), that the $c_n^{(k)}$ with fixed $k>0$ form a strictly increasing sequence in $n$. However, the sequence of the quotients $c_n^{(k)}/n$, with fixed $k\geq 0$, is decreasing as $n$ increases from 1, as the integral representation (\ref{IIIi7}) shows.  Besides, it follows from (\ref{IIIi6}) and from
\[
\lim_{n \to\infty}\left(H_n - \ln n\right) = \gamma
\]
that, for finite $k$,
\begin{equation}
\lim_{n\to\infty}\frac{c_n^{(k)}}{(\ln n)^k} = \frac{1}{k!}\,.   \label{IIIi11}
\end{equation}

\noindent{\bf Remark 4.}
The integral representation (\ref{IIIi7}) allows to establish the relation between the Roman harmonic numbers and the Coffey's probability distribution moments defined in (\ref{C1}) and (\ref{C2}). The change of variable $\delta=-(1/2)\ln(1-x)$ in these equations gives
\begin{equation}
\bar{d}^p = \frac{(-1)^p}{2^p}\,\left(D_\varepsilon -1\right)\int_0^1 dx\,x^{D_\varepsilon-2}\,\left[\ln(1-x)\right]^p\,,  \label{C3}
\end{equation}
from which, by comparison with (\ref{IIIi7}), one obtains the relation
\begin{equation}
\bar{d}^p = \frac{p!}{2^p}\,c_{D_\varepsilon -1}^{(p)}\,.   \label{C4}
\end{equation}
In view of this, expressions for the $c_n^{(k)}$ can be obtained from those of $\bar{d}^p$ given in Ref.~\cite{cof1}. For instance, Eq. (\ref{II18}) stems from Eq.~(2a) of \cite{cof1}. An interesting expansion of $c_n^{(k)}$ in series of Stirling numbers of the second kind,
\begin{equation}
c_n^{(k)}=\frac{(n-1)!}{k!}\sum_{j=n-1}^\infty\frac{(k+j)!}{j!\,n^{k+j}}\,S(j,n-1)\,,   \label{C5}
\end{equation}
results from Eq.~(5) of \cite{cof1}.

\subsection{Binomial transform}

There exist different definitions of binomial transform. Here we adopt that due to Knuth \cite {knut} and recalled in Wikipedia. Accordingly, a sequence $\{b_n\}$,  ($n=0, 1, 2, \dots$), is the binomial transform of another $\{a_n\}$ if
\begin{equation}
b_n=\sum_{l=0}^n\,(-1)^l\,{n \choose l}\,a_l\,. \label{IIIii1}
\end{equation}
The inverse transform is the same binomial one, that is,
\begin{equation}
a_n=\sum_{j=0}^n (-1)^j\,{n \choose j}\,b_j\,. \label{IIIii2}
\end{equation}
The ordinary generating functions of those sequences,
\begin{equation}
f(z)=\sum_{n=0}^\infty a_n\,z^n \quad \mbox{and} \quad g(z)=\sum_{n=0}^\infty b_n\,z^n\,,   \label{IIIii3}
\end{equation}
are related by the Euler transformation
\begin{equation}
g(z)=\frac{1}{1-z}\,f\left(\frac{z}{z-1}\right)\,,  \label{IIIii4}
\end{equation}
whereas their exponential generating functions
\begin{equation}
F(z)=\sum_{n=0}^\infty a_n\,\frac{z^n}{n!} \quad \mbox{and} \quad G(z)=\sum_{n=0}^\infty b_n\,\frac{z^n}{n!}   \label{IIIii5}
\end{equation}
are related by
\begin{equation}
G(z)=e^z\,F(-z)\,.  \label{IIIii6}
\end{equation}

\begin{lem}\hspace{-4pt}{\bf .}
The sequence of  numbers $\{c_{n+1}^{(k)}/(n+1)\}$, (n=0, 1, 2, \ldots) is the binomial transform of the sequence $\{1/(n+1)^{k+1}\}$, (n=0, 1, 2, \ldots).
\end{lem}

\noindent \textsc{Proof.}
Define
\begin{equation}
b_n=\frac{c_{n+1}^{(k)}}{n+1}\, \quad {\rm and} \quad a_n=\frac{1}{(n+1)^{k+1}}\,,  \qquad  n=0, 1, 2, \ldots\,.   \label{IIIii7}
\end{equation}
Then, Eq.~(\ref{II18}) becomes
\begin{equation}
b_{n-1}=\sum_{l=0}^{n-1}(-1)^l\,{n-1 \choose l}\,a_l\,,    \label{IIIii8}
\end{equation}
to be compared with (\ref{IIIii1}).        \hspace{212pt} $\square$

\begin{lem}\hspace{-4pt}{\bf .}
The sequence of numbers  $\{0, - c_{n}^{(k-1)}/n\}$, (n=1, 2, 3, \ldots) is the binomial transform of the sequence $\{0,H_n^{(k)}\}$, (n=1, 2, 3, \ldots).
\end{lem}

\noindent \textsc{Proof.}
Define now
\begin{equation}
a_0=0\,, \qquad a_n=H_n^{(k)}\,.   \label{IIIii9}
\end{equation}
Then, from (\ref{IIIii1}),
\begin{eqnarray}
b_0=0\,, \qquad b_n&=&\sum_{l=1}^n(-1)^l{n \choose l}\sum_{j=1}^l\frac{1}{j^k}  \nonumber  \\
&=&\sum_{j=1}^n\frac{1}{j^k}\sum_{l=j}^n(-1)^l{n \choose l}  \nonumber \\
&=&-\,\sum_{j=1}^n\frac{1}{j^k}\sum_{l=0}^{j-1}(-1)^l{n \choose l}  \label{IIIii10}
\end{eqnarray}
The fact that \cite[Eq.~(213)]{verm}
\begin{equation}
\sum_{l=0}^{j-1}(-1)^l{n \choose l} = (-1)^{j-1}{n-1 \choose j-1} =\frac{j}{n}(-1)^{j-1}{n \choose j}  \nonumber
\end{equation}
allows to complete the proof.   \hspace{199pt} $\square$

Lemma 4 may be seen as a mere consequence of Lemma 3. The Problem 6 at the end of Ch.~1 of Ref. \cite{rio2} illustrates the fact.

\subsection{Generating relations}

Since the $c_n^{(k)}$ are labeled by two indices, three kinds of generating functions may be considered: a sum over the order $k$, a sum over the degree $n$, or a double sum over both indices. Examples of generating functions of the three types are presented in the following

\begin{prop}\hspace{-4pt}{\bf .}
The Roman harmonic numbers $c_n^{(k)}$ of positive degree, $n>0$, have the following generating functions
\begin{description}
\item[{\rm (a)}]
\begin{equation}
\frac{n!}{(1-z)_n}=\sum_{k=0}^\infty c_n^{(k)}\,z^k\,, \qquad |z|<1\,.  \label{IIIiii1}
\end{equation}
\item[{\rm (b)}]
\begin{equation}
z\,e^z\,\ _{k+1}\! F_{k+1} \left( \!\! \left.\begin{array}{l}1,1,\ldots,1 \\ 2,2,\ldots,2\end{array}\right|-z\right) =
\sum_{n=1}^\infty \frac{1}{n!}\,c_n^{(k)}\,z^n\,.  \label{IIIiii2}
\end{equation}
\item[{\rm (c)}]
\begin{equation}
-\,{\rm Li}_{k+1}\left(\frac{z}{z-1}\right)=\sum_{n=1}^\infty \frac{1}{n}\,c_n^{(k)}\,z^n\,, \qquad |z|<1\,.  \label{IIIiii3}
\end{equation}
\item[{\rm (d)}]
\begin{equation}
\frac{1}{z-1}\,{\rm Li}_{k}\left(\frac{z}{z-1}\right)=\sum_{n=1}^\infty \,c_n^{(k)}\,z^n\,, \qquad |z|<1\,.   \label{IIIiii4}
\end{equation}
\item[{\rm (e)}]
\begin{equation}
\ _2\! F_1(1,1;1-z;t)=\sum_{n=0}^\infty\,\sum_{k=0}^\infty c_n^{(k)}\,z^k\,t^n\,, \qquad |z|<1\,, \quad |t|<1\,.  \label{IIIiii5}
\end{equation}
\item[{\rm (f)}]
\begin{equation}
\ _1\! F_1(1;1-z;t)=\sum_{n=0}^\infty\,\sum_{k=0}^\infty \frac{1}{n!}\,c_n^{(k)}\,z^k\,t^n\,, \qquad |z|<1\,.  \label{IIIiii6}
\end{equation}
\item[{\rm (g)}]
\begin{equation}
\ _0\! F_1(1-z;t)=\sum_{n=0}^\infty\,\sum_{k=0}^\infty \frac{1}{(n!)^2}\,c_n^{(k)}\,z^k\,t^n\,, \qquad |z|<1\,.  \label{IIIiii7}
\end{equation}
\end{description}
\end{prop}

\noindent \textsc{Proof.}
(a) becomes trivial, bearing in mind the second of (\ref{II5}) and (\ref{II13}). (b) is a consequence of the relation (\ref{IIIii6}), in view of Lemma 3 and the exponential generating function of the sequence  $\{1/(n+1)^{k+1}\}$, $(n=0, 1, 2, \ldots)$. (c) follows immediately from  Lemma 3, the definition of the polylogarithm function Li, and the Euler transformation (\ref{IIIii4}). Or, alternatively, from Lemma 4, the evident relation (assuming $H_0^{(k)}=0$)
\[
\sum_{n=0}^\infty H_n^{(k)}\,z^n=\frac{1}{1-z}\,{\rm Li}_k(z)\,, \qquad k\geq 1\,,
\]
and the Euler transformation.
(d) can be obtained from (\ref{II11}) divided by $n$, multiplied by $z^n$, and summed for $n$ from $1$ to $\infty$. Or by derivation of both sides of (\ref{IIIiii3}) with respect to $z$. (e), (f), and (g) result from using (\ref{IIIiii1}) in the corresponding series expansions of the hypergeometric functions.         \hspace{125pt} $\square$

\noindent{\bf Remark 5.}
The left-hand side of (\ref{IIIiii1}) can be written in terms of the beta function \cite[Eq.~(A.1)]{coff} to give
\begin{equation}
n\,B(1-z,n)=\sum_{k=0}^\infty c_n^{(k)}\,z^k\,, \qquad |z|<1\,.  \label{IIIiii8}
\end{equation}
Then,
\begin{equation}
c_n^{(k)}=(-1)^k\frac{n}{k!}\left[\frac{\partial^k}{\partial z^k} B(z,n)\right]_{z=1}\,, \label{IIIiii9}
\end{equation}
in accordance with (\ref{IIIi7}), in view of \cite[Eq.~2.6.9.5]{prud}.

\noindent{\bf Remark 6.}
By using in the right-hand side of (\ref{IIIiii2}) the integral representation (\ref{IIIi7}) one obtains easily
\begin{equation}
\ _{k+1}\! F_{k+1} \left( \!\! \left.\begin{array}{l}1,1,\ldots,1 \\ 2,2,\ldots,2\end{array}\right|z\right) =
\frac{1}{k!}\int_0^1dx\,e^{zx}\,(-\ln x)^k\,.      \label{IIIiii10}
\end{equation}

\noindent{\bf Remark 7.}
Two series expansions of the polylogarithm function, convergent in a half-plane, namely
\begin{eqnarray}
{\rm Li}_{k}(t) =  -\,\sum_{n=1}^\infty \frac{1}{n}\,c_n^{(k-1)}\,\left(\frac{t}{t-1}\right)^n\,, \qquad &\Re\, t<1/2\,,& \quad k\geq 1\,,   \label{IIIiii11}  \\
{\rm Li}_{k}(t) = \frac{1}{t-1}\,\sum_{n=1}^\infty \,c_n^{(k)}\,\left(\frac{t}{t-1}\right)^n\,, \qquad &\Re\, t<1/2\,,&  \quad k\geq 0\,,   \label{IIIiii12}
\end{eqnarray}
can be immediately obtained from (\ref{IIIiii3}) and (\ref{IIIiii4}).

\noindent{\bf Remark 8.}
The relation (\ref{IIIiii7}) suggests a representation of the Bessel function $J_\nu(z)$ which may be useful for computing the successive derivatives with respect to its order $\nu$, a problem which has deserved interest recently \cite{bry2}. Such representation reads
\begin{equation}
J_\nu(z) = \frac{(z/2)^\nu}{\Gamma(1+\nu)}\sum_{k=0}^\infty \mathcal{F}_k(-z^2/4)\,(-\nu)^k\,,\qquad  |\nu|<1\,, \label{IIIiii13}
\end{equation}
having defined
\begin{equation}
\mathcal{F}_k(z) = \sum_{n=0}^\infty\frac{c_n^{(k)}}{(n!)^2}\,z^n\,.   \label{IIIiii14}
\end{equation}
The functions $\{\mathcal{F}_k\}$ ($k=0, 1, \ldots$) obey the differential-difference equation
\begin{equation}
\left(\frac{d}{dz}\,z\,\frac{d}{dz}-1\right)\mathcal{F}_k(z)=\frac{d}{dz}\mathcal{F}_{k-1}(z)\,.  \label{IIIiii15}
\end{equation}
Obviously,
\begin{equation}
\mathcal{F}_0(z) = \,_0\!F_1(1;z) = I_0\!\left(2\sqrt{z}\right)  \label{IIIiii16}
\end{equation}
and, in view of (\ref{II19}),
\begin{equation}
\lim_{k\to\infty} \mathcal{F}_k(z) =z\,_0\!F_1(2;z) = \sqrt{z}\,I_1\!\left(2\sqrt{z}\right)\,,  \label{IIIiii17}
\end{equation}
where $I$ stands for the modified Bessel function.
Besides, from the relation \cite[Vol.~3, Eq.~7.13.1.6]{prud}
\begin{equation}
\,_0\!F_1\left(\frac{3}{2};t\right) = \frac{\sinh \left(2\sqrt{t}\right)}{2\sqrt{t}}  \label{IIIiii18}
\end{equation}
and Eq.~(\ref{IIIiii7}) with $z$ replaced by $-1/2$, one obtains
\begin{equation}
\frac{\sinh \left(2\sqrt{t}\right)}{2\sqrt{t}}  = \sum_{k=0}^\infty \frac{(-1)^k}{2^k}\sum_{n=0}^\infty\frac{c_n^{(k)}}{(n!)^2}\,t^n\,,   \label{IIIiii19}
\end{equation}
which allows one to write the sum rule for the family of functions $\{\mathcal{F}_k\}$
\begin{equation}
\sum_{k=0}^\infty \frac{(-1)^k}{2^k}\,\mathcal{F}_k(z) = \frac{\sinh \left(2\sqrt{z}\right)}{2\sqrt{z}}\,.   \label{IIIiii20}
\end{equation}

\subsection{Sum rules}

Sum rules for the Roman harmonic numbers of negative degree can be immediately written by translation of known sum rules for the Stirling numbers of the first kind. For instance, the sum rules \cite[Eqs.~26.8.28 and 26.8.29]{dlmf}
\begin{eqnarray}
\sum_{k=1}^\infty\,s(n,k) = 0\,, \qquad &{\rm for}& \quad n\geq 2\,,   \label{IIIiv1a}  \\
\sum_{k=1}^\infty\,(-1)^{n-k}\,s(n,k) = n!\,, \qquad &{\rm for}& \quad n\geq 0\,,   \label{IIIiv1b}
\end{eqnarray}
are translated into
\begin{eqnarray}
\sum_{k=1}^\infty\,(-1)^{k}\,c_{-n}^{(k)}  =  0\,, \qquad  &{\rm for}& \quad n\geq 2\,,   \label{IIIiv2a}  \\
\sum_{k=1}^\infty\,c_{-n}^{(k)}  =  -\,n\,, \qquad  &{\rm for}& \quad n\geq 0\,.  \label{IIIiv2b}
\end{eqnarray}
Notice that, since $s(n,k)=c_{-n}^{(k)}=0$ for $k>n\geq 0$, the sums in (\ref{IIIiv1a}) to (\ref{IIIiv2b}) are finite.
In the case of positive degree, instead, the property (c) of Proposition 2 makes evident that the infinite series
\[
\sum_{k=0}^\infty\,c_n^{(k)}  \qquad {\rm and} \qquad \sum_{k=0}^\infty\,(-1)^k\,c_n^{(k)}
\]
are respectively divergent and non-convergent. Nevertheless, the successive partial sums of the second of these series oscillate in the interval $[1/(n+1)-n/2, 1/(n+1)+n/2]$ and the series, although non-convergent, may be summed in the sense of Abel \cite{hard} to give
\begin{equation}
\sum_{k=0}^\infty\,(-1)^k\,c_n^{(k)}=1/(n+1)\,, \qquad {\rm (A)}\,. \label{IIIiv3}
\end{equation}
The same result is obtained, formally, by considering $z\to -1^+$ in (\ref{IIIiii1}) or by using the integral representation (\ref{IIIi7}).

A sum rule with binomial coefficients follows immediately from Lemma 3. It reads
\begin{equation}
\sum_{j=1}^n\,(-1)^{j-1}\,{n\choose j}\,c_j^{(k)} = 1/n^k\,.        \label{IIIiv4}
\end{equation}

The generating relations shown above allow to write possibly useful sum rules. For instance, by taking respectively $z=1/2$ and $z=-1/2$ in (\ref{IIIiii1}) one has
\begin{equation}
\sum_{k=0}^\infty\,\frac{c_n^{(k)}}{2^k}=\frac{2^{2n}\,(n!)^2}{(2n)!}\,, \qquad
\sum_{k=0}^\infty\,(-1)^k\,\frac{c_n^{(k)}}{2^k}=\frac{2^{2n}\,(n!)^2}{(2n+1)!}\,, \qquad n\geq 0\,, \label{IIIiv5}
\end{equation}
from which
\begin{equation}
\sum_{k=0}^\infty\,\frac{c_n^{(2k)}}{2^{2k}}=\frac{2^{2n}\,n!\,(n+1)!}{(2n+1)!}\,, \qquad
\sum_{k=0}^\infty\,\frac{c_n^{(2k+1)}}{2^{2k+1}}=\frac{2^{2n}\,n\,(n!)^2}{(2n+1)!}\,, \qquad n\geq 0\,, \label{IIIiv6}
\end{equation}
are immediately obtained. A similar procedure, with different values of $z$, can be used to obtain a variety of sum rules from (\ref{IIIiii3}) and (\ref{IIIiii4}). For example,
\begin{equation}
\sum_{n=0}^\infty\,\frac{c_n^{(k)}}{2^n\,n}=-{\rm Li}_{k+1}(-1)\,, \qquad\qquad
\sum_{n=0}^\infty\,\frac{c_n^{(k)}}{2^n}=-2\,{\rm Li}_{k}(-1)\,. \label{IIIiv7}
\end{equation}
In the case of the label of the polylogarithm function being even, one may benefit from the relation \cite[Appendix II.5]{prud}
\begin{equation}
{\rm Li}_{2k}(-1) = (-1)^k\,\frac{2^{2k-1}-1}{(2k)!}\,\pi^{2k}\,B_{2k}\,,  \label{IIIiv8}
\end{equation}
where $B_{2k}$ represents a Bernoulli number.

\begin{lem}\hspace{-4pt}{\bf .}
The Stirling numbers of the first kind obey the sum rule
\begin{equation}
\sum_{j=0}^k\,(-1)^j\,s(n+1,k+1-j)\,s(-n,j)=(-1)^n\,\delta_{k,0}\,, \qquad n\geq 0\,.   \label{IIIiv9}
\end{equation}
\end{lem}

\noindent \textsc{Proof.}
Consider the evident identity
\begin{equation}
\frac{\Gamma(x+1)}{\Gamma(x-n)}\,\frac{\Gamma(-x+1)}{\Gamma(-x+1+n)} = (-1)^n\,x\,,  \qquad n\geq 0\,, \label{IIIiv10}
\end{equation}
expand the two fractions on the left-hand side according to (\ref{II4}), and compare coefficients of $x^{k+1}$ in the two sides. \hspace{188pt} $\square$

\begin{prop}\hspace{-4pt}{\bf .}
The Roman harmonic numbers obey the following sum rules
\begin{description}
\item[{\rm (a)}]
\begin{equation}
\sum_{j=0}^k\,s(n+1,k+1-j)\,c_n^{(j)}=(-1)^n\,\delta_{k,0}\,n!\,, \qquad n\geq 0\,.    \label{IIIiv11}
\end{equation}
\item[{\rm (b)}]
\begin{equation}
\sum_{j=0}^k (-1)^j\,c_{-n-1}^{(k+1-j)}\,c_n^{(j)}=-\,\delta_{k,0}\,, \qquad n\geq 0\,.  \label{IIIiv12}
\end{equation}
\item[{\rm (c)}]
\begin{equation}
\sum_{j=0}^k \,\frac{B_{n-k+j}^{(n+1)}}{(n-k+j)!}\,\frac{c_n^{(j)}}{(k-j)!}=(-1)^n\,\delta_{k,0}\, \qquad n\geq 0\,.
    \label{IIIiv13}
\end{equation}
\end{description}
\end{prop}

\noindent \textsc{Proof.}
Substitution of (\ref{II13}) in (\ref{IIIiv9}) proves immediately (a) and (b). (c) is obtained from (a) by using the relation between Stirling and generalized Bernoulli numbers \cite[Eq.~(2.28)]{coff}
\begin{equation}
s(m,n)= {m-1 \choose n-1}\, B_{m-n}^{(m)}\,.
\end{equation}
 \hspace{334pt}   $\square$

\noindent{\bf Remark 9.}
Interesting particular cases of the above sum rules are
\begin{eqnarray}
\sum_{j=0}^n\,s(n+1,n+1-j)\,c_n^{(j)}&=&\delta_{n,0}\,, \qquad n\geq 0\,,    \label{IIIiv14}   \\
\sum_{j=0}^n (-1)^j\,c_{-n-1}^{(n+1-j)}\,c_n^{(j)}&=&-\,\delta_{n,0}\,, \hspace{-11pt}\qquad n\geq 0\,,
\label{IIIiv15}  \\
\sum_{j=0}^n \,\frac{B_{j}^{(n+1)}\,c_n^{(j)}}{j!\,(n-j)!}&=&\delta_{n,0}\,, \qquad n\geq 0\,.   \label{IIIiv16}
\end{eqnarray}

\noindent{\bf Remark 10.}
The ``orthogonality" relations (\ref{IIIiv11}) to (\ref{IIIiv13}) allow to write inverse relations \cite[Ch.~2 and 3]{rio2} between two sequences $\{a_m\}$ and $\{b_m\}$. In fact, for arbitrary $ n\geq 0$, we may state that
\begin{equation}
{\rm if}  \quad b_m=\sum_{l=0}^m s(n+1,l+1)\,a_{m-l}\,,  \quad {\rm then} \quad a_m=\frac{(-1)^n}{n!}\sum_{j=0}^m c_n^{(j)}\,b_{m-j}\,,     \label{IIIiv17}
\end{equation}
or, in terms of generalized Bernoulli numbers,
\begin{equation}
{\rm if}  \quad b_m=\sum_{l=0}^m {n \choose l}B_{n-l}^{(n+1)}\,\,a_{m-l}\,,  \quad {\rm then} \quad a_m=\frac{(-1)^n}{n!}\sum_{j=0}^m c_n^{(j)}\,b_{m-j}\,.     \nonumber
\end{equation}
Analogously,
\begin{equation}
{\rm if}  \quad b_m=\sum_{l=0}^m (-1)^{l+1}\,c_{-n-1}^{(l+1)}\,a_{m-l}\,,  \quad {\rm then} \quad a_m=\sum_{j=0}^m c_n^{(j)}\,b_{m-j}\,.     \label{IIIiv18}
\end{equation}

\subsection{Derivatives of the Pochhammer symbol and of its reciprocal}

In a previous paper \cite{grey}, we have shown that having expressions for the derivatives, to any order, of the Pochhammer and reciprocal Pochhammer  symbols with respect to their arguments facilitates the evaluation of the so called $\varepsilon$-expansion of functions of the hypergeometric class related to Feynman integrals in high-energy physics. Equations (\ref{II5}) and (\ref{II13}) suggest the possibility of
writing such derivatives in terms of Roman harmonic numbers. In the case of the Pochhammer symbol, from
\begin{equation}
(x)_n = -\,(n-1)! \sum_{l=0}^n c_{-n}^{(l)}\,x^l\,,  \qquad  n\geq 1\,,    \label{IIIv0}
\end{equation}
one obtains immediately
\begin{equation}
\frac{d^k}{dx^k}(x)_n = -\,(n-1)!\sum_{j=0}^{n-k}(j+1)_k\,c_{-n}^{(k+j)}\,x^j\,, \qquad  n\geq 1\,,   \label{IIIv1}
\end{equation}
which, by making use of (\ref{II15}) and (\ref{IIIv0}), can be written in the form
\begin{equation}
\frac{d^k}{dx^k}(x)_n = n!\,k!\sum_{j=0}^{n-k}\frac{c_{-n+j}^{(k)}}{-n+j}\,\frac{(x)_j}{j!}\,, \qquad  n\geq 1\,,   \label{IIIv2}
\end{equation}
Especially interesting for the computation of the mentioned Feynman integrals is the case of the argument of the derivatives being an integer. Obviously,
\begin{eqnarray}
\left. \frac{d^k}{dx^k}(x)_n\right|_{x=0} &=& -\,(n-1)!\,k!\,c_{-n}^{(k)}\,, \qquad  n\geq 1\,,   \label{IIIv3}  \\
\left. \frac{d^k}{dx^k}(x)_n\right|_{x=1} &=& n!\,k!\sum_{j=0}^{n-k}\frac{c_{-n+j}^{(k)}}{-n+j}\,, \qquad  n\geq 1\,.  \label{IIIv4}
\end{eqnarray}
In view of (\ref{II16}), the last equation can be written in the form
\begin{equation}
\left. \frac{d^k}{dx^k}(x)_n\right|_{x=1} = -\,n!\,k!\,c_{-n-1}^{(k+1)}\,. \qquad  n\geq 1\,.  \label{IIIv5}
\end{equation}
Expressions for other integer values of the argument, alternative to those stemming from (\ref{IIIv1}) or (\ref{IIIv2}), will be given below.

In what concerns the reciprocal Pochhammer symbol, we obtain from (\ref{II5}) and (\ref{II13}) the series expansion
\begin{equation}
\frac{1}{(x)_n} = \frac{1}{(n-1!}\,\frac{1}{x}\sum_{j=0}^\infty(-1)^j\,c_{n-1}^{(j)}\,x^j\,,  \qquad 0<|x|<1\,,  \label{IIIv6}
\end{equation}
and, by repeated derivation with respect to $x$,
\begin{equation}
\frac{d^k}{dx^k} \frac{1}{(x)_n} = \frac{(-1)^k}{(n-1)!}\left(\frac{k!}{x^{k+1}}+\sum_{l=0}^\infty(-1)^{l+1}\,(l+1)_k\,c_{n-1}^{(k+l+1)}\,x^l\right),  \quad 0<|x|<1\,.  \label{IIIv7}
\end{equation}
For values of $x$ out of the open disk, i. e., for $|x|\geq 1$, use can be made of the partial fraction decomposition \cite[Vol 1, Eq.~4.2.2.45]{prud}
\begin{equation}
\frac{1}{(x)_n} = \sum_{j=0}^{n-1}\frac{(-1)^j}{j!\,(n-1-j)!}\,\frac{1}{x+j}\,, \qquad n\geq 1\,.  \label{IIIv8}
\end{equation}
This representation of the function $1/(x)_n$ is valid in the whole complex $x$-plane with exception of the points $x=0, -1, \ldots, -(n-1)$, where it presents poles of first order. By deriving repeatedly with respect to $x$, one obtains
\begin{equation}
\frac{d^k}{dx^k} \frac{1}{(x)_n} = \frac{(-1)^k\,k!}{(n-1)!}\,\sum_{j=0}^{n-1}(-1)^j\,{n-1 \choose j} \frac{1}{(x+j)^{k+1}}\,,  \label{IIIv9}
\end{equation}
an expression which shows that, for every non-negative integer $k$, the sequences
\[
\left\{n!\,\frac{d^k}{dx^k}\frac{1}{(x)_{n+1}}\right\} \quad {\rm and} \quad \left\{\frac{(-1)^k\,k!}{(x+n)^{k+1}}\right\}\,,
\quad n=0, 1, 2, \ldots\,,
\]
are binomial transform of each other. It also shows that
\begin{equation}
\frac{d^k}{dx^k} \frac{1}{(x)_n} = \frac{(-1)^{n-k+1}\,k!}{(n-1)!}\,\Delta^{n-1}\frac{1}{x^{k+1}}\,, \quad n>0\,, \label{IIIv10}
\end{equation}
where $\Delta$ represents the difference operator \cite[Sec.~1.6]{comt} \cite[Sec.~6.2]{rio2}, defined by
\begin{equation}
\Delta^n f(x)=\sum_{j=0}^n\,(-1)^{n-j}\,{n \choose j}\,f(x+j)\,.  \label{IIIv11}
\end{equation}

A useful test for the expressions given for the derivatives of the Pochhammer and reciprocal Pochhammer symbols can be immediately obtained from the identity
\begin{equation}
(x)_n\,\frac{1}{(x)_n} = 1\,.
\end{equation}
It reads
\begin{equation}
\sum_{j=0}^k {k \choose j} \left(\frac{d^{k-j}}{dx^{k-j}}(x)_n\right)\left(\frac{d^j}{dx^j}\frac{1}{(x)_n}\right) = \delta_{k,0}\,.
\end{equation}

We have already mentioned our interest in having expressions of the derivatives of the Pochhammer and reciprocal Pochhammer symbols for integer values of the variable. From (\ref{IIIv9}), we see that
\begin{equation}
\left.\frac{d^k}{dx^k} \frac{1}{(x)_n} \right|_{x=1} = \frac{(-1)^k\,k!}{n!}\,c_n^{(k)}\,,  \label{IIIv12}
\end{equation}
which, together with (\ref{IIIv3}) and (\ref{IIIv5}), provides what we need for our purpose. Bearing in mind these equations, and the identity
\begin{equation}
(x)_n=(x-N+1)_{n+N-1}\,\frac{1}{(x-N+1)_{N-1}} \quad {\rm for} \quad N>0\,,  \label{IIIv13}
\end{equation}
it is easy to obtain, for $\quad N>0\,,$
\begin{eqnarray}
\left.\frac{d^k}{dx^k} (x)_n \right|_{x=N} &=& k!\,(N)_n\,\sum_{j=0}^k (-1)^{j+1}\,c_{-n-N}^{(k+1-j)}\,c_{N-1}^{(j)}\,,  \label{IIIv14}  \\
\left.\frac{d^k}{dx^k} \frac{1}{(x)_n} \right|_{x=N} &=& \frac{k!}{(N)_n}\,\sum_{j=0}^k (-1)^{j+1}\,c_{-N}^{(k+1-j)}\,c_{n+N-1}^{(j)}\,,  \label{IIIv15}
\end{eqnarray}
Analogously, since
\begin{equation}
(x)_n=(-1)^n\,(-x-n+1)_n\,,   \label{IIIv16}
\end{equation}
we can write, for $\quad N\geq n>0\,,$
\begin{equation}
\left.\frac{d^k}{dx^k} (x)_n \right|_{x=-N} = (-1)^{n-k}\,k!\,(N-n+1)_n\,\sum_{j=0}^k (-1)^{j+1}\,c_{-N-1}^{(k+1-j)}\,c_{N-n}^{(j)}\,,
\label{IIIv17}
\end{equation}
\begin{equation}
\left.\frac{d^k}{dx^k} \frac{1}{(x)_n} \right|_{x=-N} = (-1)^{n-k}\,\frac{k!}{(N-n+1)_n}\,\sum_{j=0}^k (-1)^{j+1}\,c_{-N+n-1}^{(k+1-j)}\,c_{N}^{(j)}\,. \label{IIIv18}
\end{equation}
In the case of being $\quad n>N>0\,,\quad$ we may write
\begin{equation}
(x)_n=(-1)^{N}\,(-x-N+1)_{N}\,(x+N)_{n-N}\,,  \label{IIIv19}
\end{equation}
and therefore
\begin{equation}
\left.\frac{d^k}{dx^k} (x)_n \right|_{x=-N} = (-1)^{N-k}\,k!\,N!\,(n-N-1)!\,\sum_{j=0}^k (-1)^{j}\,c_{-N-1}^{(k-j+1)}\,c_{-n+N}^{(j)}\,.
\label{IIIv20}
\end{equation}

Alternative expressions for the derivatives of the reciprocal Pochhammer symbol can be obtained with a quite different strategy. The Fa\`{a} di Bruno's formula \cite[Secs.~3.3 and 3.4]{comt}
\begin{equation}
\frac{d^k}{dx^k}f(g(x)) = \sum_{j=1}^k \left.\frac{d^j}{dt^j}f(t)\right|_{t=g(x)}\,\mathbf{B}_{k,j}\left(\frac{d}{dx}g(x),\,\frac{d^2}{dx^2}g(x), \ldots, \frac{d^{k-j+1}}{dx^{k-j+1}}g(x)\right),   \label{IIIv21}
\end{equation}
where $\mathbf{B}_{k,j}$ stands for the partial exponential Bell polynomials, allows one to write
\begin{equation}
\frac{d^k}{dx^k}\frac{1}{(x)_n}= \sum_{j=1}^k \frac{(-1)^j\,j!}{\left((x)_n\right)^{j+1}}\,\mathbf{B}_{k,j}\left(\frac{d}{dx}(x)_n,\,\frac{d^2}{dx^2}(x)_n, \ldots, \frac{d^{k-j+1}}{dx^{k-j+1}}(x)_n\right),   \label{IIIv22}
\end{equation}
which gives the the derivatives of the reciprocal Pochhammer symbol in terms of those of the Pochhammer symbol. Particularizing both sides of this relation for $x=1$, one has
\begin{eqnarray}
c_n^{(k)} &=& \frac{(-1)^k}{k!}\sum_{j=1}^k\frac{(-1)^{j}\,j!}{(n!)^{j}}   \nonumber  \\
 & & \hspace{-25pt}\times\,\mathbf{B}_{k,j}\left(-n!1!c_{-n-1}^{(2)}, -n!2!c_{-n-1}^{(3)},\ldots, -n!(k-j+1)!c_{-n-1}^{(k-j+2)}\right),   \label{IIIv23}
\end{eqnarray}
an expression of the positive degree Roman harmonic numbers in terms of the negative degree ones. By using unsigned Stirling numbers and
the partial ordinary Bell polynomials $\hat{\mathbf{B}}_{k,j}$, the last relation becomes
\begin{equation}
c_n^{(k)} = \sum_{j=1}^k\frac{(-1)^{k-j}}{(n!)^{j}}\,\hat{\mathbf{B}}_{k,j}\left(|s(n+1,2)|, |s(n+1,3|, \ldots , |s(n+1,k-j+2)|\right).  \label{IIIv24}
\end{equation}

\section{Final comments}

In the preceding sections we have shown a collection of properties of the Roman harmonic numbers. This paper does not pretend to be an exhaustive revision of the theme. Quite the contrary, our purpose was to stimulate the potential readers to discover new properties and applications of those interesting numbers, forgotten for a long time.

Since they are rational numbers, it is possible to evaluate and store the Roman harmonic numbers in an exact way. This fact makes them especially suitable for use in computational tasks.

\section*{Acknowledgements}

The author is very indebted to an anonymous reviewer for pointing out the connection of the Coffey's probability distribution moments with the Roman harmonic numbers, and contributing, in this way, to make this paper more interesting. This work has been supported by Departamento de Ciencia, Tecnolog\'{\i}a y Universidad del Gobierno de Arag\'on (Project 226223/1) and Ministerio de Ciencia e Innovaci\'on (Project MTM2015-64166).

\end{document}